\documentclass{elsarticle}

\usepackage{amsmath}
\usepackage{bm}
\usepackage{tikz}
\usepackage{pgfplots}

\usepackage{hyperref}

\newtheorem{thm}{Theorem}

\newproof{pf}{Proof}

\journal{arXiv.org}

\begin{document}

\begin{frontmatter}

\title{Numerical solving the boundary value problem for fractional powers of elliptic operators}

\author[mymainaddress,mysecondaryaddress]{Petr N. Vabishchevich}
\ead{vabishchevich@gmail.com}

\address[mymainaddress]{Nuclear Safety Institute, Russian Academy of Sciences, 52, B. Tulskaya, Moscow, Russia}
\address[mysecondaryaddress]{North-Eastern Federal University, 58, Belinskogo, Yakutsk, Russia}

\begin{abstract}
A boundary value problem for a fractional power of the second-order elliptic operator is considered.
It is solved numerically using a time-dependent problem for a pseudo-parabolic equation.
For the auxiliary Cauchy problem, the standard two-level schemes with weights are applied.
Stability conditions are obtained for the fully discrete schemes under the consideration.
The numerical results are presented for a model two-dimensional
boundary value problem wit a fractional power of an elliptic operator.
The dependence of accuracy on grids in time and in space is studied.
\end{abstract}


\end{frontmatter}

\section{Introduction}

Nowadays, non-local applied mathematical models based on the use of fractional derivatives in time and space are actively discussed \cite{eringen2002nonlocal,kilbas2006theory,baleanu2012fractional}. 
An interesting example is a boundary value problem for a fractional power of an elliptic operator. 
For example, in a bounded domain $\Omega$, we search the solution of the problem
\[
 (- \triangle )^\alpha  u = f(\bm x), 
 \quad \bm x \in \Omega , 
\] 
\[
  u(\bm x) = 0,
 \quad \bm x \in \partial \Omega , 
\] 
where $0 < \alpha < 1$.

Different approaches are employed to solve numerically such boundary value problems.
The simplest variant is associated with the explicit construction of  the solution using the known 
eigenvalues ​​and eigenfunctions of the elliptic operator with 
diagonalization of the corresponding matrix 
\cite{ilic2005numerical,ilic2006numerical,yang2011novel}. 
Unfortunately, this approach demonstrates too high computational complexity for multidimensional problems.

Another approach is based on the representation of an elliptic operator power
in the form of a contour integral (the Dunford-Cauchy representation)
and the application of appropriate quadrature formulas with nodes of integration on the complex
plane. The approximate operator is written as the sum of resolvents \cite{gavrilyuk2004data,gavrilyuk2005data}
that provides the exponential convergence of quadrature approximations.
The paper \cite{bonito2013numerical} uses quadrature
formulas with nodes on the real axis, which are constructed
using the corresponding integral representation for the operator power 
\cite{krasnoselskii1976integral}.
In this case, for the inverse operator of the problem, we have an additive representation,
where each operator term  is the normal inverse elliptic operator.

It is also necessary to note the possibility of finding the solution for the fractional
power of an elliptic operator as the solution of an
elliptic boundary value problem of higher dimension, i.e., we introduce a new variable
\cite{caffarelli2007extension,stinga2010extension}.
In the paper \cite{nochetto2013pde}, this possibility is discussed at the differential level
along with the justification of the proposed computational algorithm.
The computational efficiency of these algorithms is not very high, and moreover,
such an extended elliptic boundary value problem is degenerate/singular and
the problem is solved in a semi-infinite cylinder.

After constructing finite difference or finite element approximations,
from the boundary value problem for the fractional power of
the elliptic operator, we arrive at the problem of 
multiplication of the fractional power of the matrix corresponding to elliptic operator
by the vector, which corresponds to the right-hand side.
For such a matrix problem, different approaches are used \cite{higham2008functions}:
Krylov subspace methods, contour integration and so on.
Special attention should be given to the methods that solve the Cauchy problem for
the corresponding evolutionary equation.
The first work in this direction for matrix problems with
$\alpha = - 0.5$  is the work \cite{allen2000numerical},
where the approximate solution is searched as  the solution of the unsteady problem
within the unit time interval.
This type of computational algorithms have obvious advantages,
related, in particular, to their relative simplicity and regularity:
the numerical solution of steady-state problems is
often based on the transition to solving pseudo-time evolutionary problems.

In the present paper, for solving the boundary
value problem for a fractional power of an elliptic operator, we use
the transition to a pseudo-parabolic equation.
The stability of the two-level scheme with weights is shown.
Numerical experiments for a model two-dimensional problem are performed
using the standard finite element approximations.
The paper is organized as follows.
The formulation of the boundary value problem for the fractional power of the elliptic operator is given in Section 2. Section 3 presents the auxiliary Cauchy problem for the pseudo-parabolic
equation. 
The central body of the work is Section 4, where time-stepping techniques are constructed and justified.
The results of numerical experiments are described in Section 5.

\section{Problem formulation} 

In a bounded polygonal domain $\Omega \subset R^m$, $m=1,2,3$ with the Lipschitz continuous boundary $\partial\Omega$, we search the solution for the problem with a fractional power of a elliptic operator.
Define the elliptic operator as
\begin{equation}\label{1}
  A u = - {\rm div}  k({\bm x}) {\rm grad} u + c({\bm x}) u
\end{equation} 
with coefficients $0 < k_1 \leq k({\bm x}) \leq k_2$, $c({\bm x}) \geq 0$.
The operator $A$ is defined on the set of functions $u({\bm x})$ that satisfy
on the boundary $\partial\Omega$ the following conditions:
\begin{equation}\label{2}
  k({\bm x}) \frac{\partial u }{\partial n } + \mu ({\bm x}) u = 0,
  \quad {\bm x} \in \partial \Omega ,
\end{equation} 
where $\mu ({\bm x}) \geq \mu_1 > 0, \  {\bm x} \in \partial \Omega$.

In the Hilbert space $H = L_2(\Omega)$, we define 
the scalar product and norm in the standard way:
\[
  <u,v> = \int_{\Omega} u({\bm x}) v({\bm x}) d{\bm x},
  \quad \|u\| = <u,u>^{1/2} .
\] 
In $H$, the operator $A$ is self-adjoint and positive definite: 
\begin{equation}\label{3}
  A = A^* \geq \delta I ,
  \quad \delta > 0 ,    
\end{equation} 
where $I$ is the identity operator in $H$.
We seek the solution of the equation with the fractional power of the operator $A$:
\begin{equation}\label{4}
  A^\alpha u = f({\bm x}),
  \quad {\bm x} \in \Omega ,   
\end{equation} 
under the restriction $0 < \alpha < 1$. 

In view of (\ref{3}),  the solution of the problem (\ref{4}) satisfies the a priori estimate
\begin{equation}\label{5}
  \|u\| \leq  \delta^{-\alpha} \|f\| ,   
\end{equation}
which is valid for all $0 < \alpha < 1$. 
For $0.5 < \alpha < 1$, the a priori estimate (\ref{5}) can be improved.
In this case, we have
\begin{equation}\label{6}
  \delta^{2\alpha} \|u\|^2 + 2 \alpha \delta^{2\alpha - 1}  <(A -\delta I) u, u>  \leq   \|f\|^2 . 
\end{equation}
In particular, for $\alpha = 0.5$, we have the identity
\begin{equation}\label{7}
  \|u\|_A^2 = \|f\|^2 . 
\end{equation}

Indeed, for the solution of equation (\ref{4}), we have
\[
 <A^{2\alpha} u,u> = <f, f>
\] 
and for $\alpha = 0.5$, we have (\ref{7}).
For $\alpha > 0.5$, in view of (\ref{3}), we apply the estimate
\[
 A^{2\alpha} = (A - \delta I + \delta I)^{2\alpha} \geq 
 \delta^{2\alpha} I + 2 \alpha \delta^{2\alpha - 1} (A - \delta I) .
\]
This leads us immediately to the a priori estimate (\ref{6}). 

\section{The Cauchy problem} 

Assume that
\[
 w(t) = \delta^{\alpha} (t (A - \delta I) + \delta I)^{-\alpha} w(0) ,
\]
then the solution of equation (\ref{4}) can be defined as
\begin{equation}\label{8}
 u = w(1),
 \quad w(0) = \delta^{-\alpha} f .  
\end{equation} 
The function $w(t)$ satisfies the evolutionary equation
\begin{equation}\label{9}
  (t D + \delta I) \frac{d w}{d t} + \alpha D w = 0 ,
\end{equation} 
where
\[
 D = A - \delta I \geq 0 .
\] 
The equation (\ref{9}) is a pseudo-parabolic equation of second order
\cite{Showalter1969,ShowalterTing1970}.

\begin{thm}\label{t-1}
The solution of the problem (\ref{8}), (\ref{9}) 
satisfies the estimate (\ref{5}) for $0 < \alpha <1$ 
or the estimate (\ref{6}) if $0.5 < \alpha <1$.
\end{thm}

\begin{pf}
Multiplying equation (\ref{9})) scalarly in $H$ by $w$, we obtain
\[
 \|w(t)\| \leq \|w(0)\| . 
\]
By (\ref{8}), this inequality results in the a priori estimate (\ref{5}) for the solution of equation 
(\ref{4}) for all $0 < \alpha < 1$.  It is a little more difficult to obtain an a priori estimate 
for the solution of the Cauchy problem for  equation (\ref{9}), which implies (\ref{6}). 
For $0.5 < \alpha < 1$, multiplying equation (\ref{9}) by
\[
  w + \frac{2\alpha -1}{\alpha } t \frac{d w}{d t} , 
\]
we get 
\[
 \frac{2\alpha -1}{\alpha } t \left <(t D + \delta I)\frac{dw}{dt}, \frac{dw}{dt}\right > + 
 \left < \left (\delta I + 2\alpha t  D \right ) \frac{dw}{dt}, w \right > +
 \alpha <Dw,w> = 0.
\]
Taking into account that 
\[
 t  \frac{dw}{dt} + \frac{1}{2} w = t^{1/2}  \frac{d}{dt} (w t^{1/2}) ,
\] 
we obtain
\[
 \frac{d}{dt} (\delta \|w\|^2 + 2 \alpha <tD w,w> ) \leq 0 .
\]
Thus we have the following a priori estimate for the solution of the Cauchy problem for  equation (\ref{9}):
\begin{equation}\label{10}
 \delta \|w(t) \|^2 + 2 \alpha <tD w(t),w(t)>  \leq \delta \|w(0) \|^2 . 
\end{equation} 
In view of (\ref{8}), from (\ref{10}), it follows the estimate (\ref{6}) for the solution of equation (\ref{4}). For $\alpha = 0.5$, we have the equality (see (\ref{7}))
\begin{equation}\label{11}
 \delta \|w(t) \|^2 + <tD w(t),w(t)> = \delta \|w(0) \|^2 
\end{equation} 
for the solution of (\ref{9}).
\end{pf}

\section{Difference scheme} 

To solve numerically the problem (\ref{8}), (\ref{9}), we employ finite element 
approximations in space \cite{ciarlet2002finite,brenner2008mathematical}. 
For (\ref{1}) and (\ref{2}), we define the bilinear form
\[
 a(u,v) = \int_{\Omega } \left ( k {\rm grad}  u \, {\rm grad}  v + c u v \right )  d {\bm x} +
 \int_{\partial \Omega } \mu u v d {\bm x} .
\] 
By (\ref{3}), we have
\[
a(u,u) \geq \delta \|u\|^2 .  
\]
Define a subspace of finite elements $V^h \subset H^p(\Omega)$,
where, for example, $p=1,2$.
The approximate solution $w \in V^h$ is defined as the solution 
of the problem
\begin{equation}\label{12}
   d \left (t\frac{ d w}{d t}, v \right ) + \delta \left < \frac{ d w}{d t}, v \right > +
  \alpha d (w,v) = 0,   
  \quad 0 < t \leq T,  
\end{equation} 
\begin{equation}\label{13}
  <w(0),  v> = <\delta^{-\alpha} f , v> ,
  \quad \forall v \in V^h ,    
\end{equation} 
where
\[
  d(u,v) \equiv a(u,v) - \delta <u,v> .
\] 

To solve numerically the problem (\ref{12}), (\ref{13}),
we use the simplest implicit two-level scheme with weights.
Let $\tau$ be a step of a uniform grid in time such that $w^n = w(t^n), \ t^n = n \tau$,
$n = 0,1, ..., N, \ N\tau = 1$.
Let us approximate equation (\ref{12}) by the implicit two-level scheme
\begin{equation}\label{14}
\begin{split}
  d \left (t^{\sigma(n)}\frac{ w^{n+1} - w^{n}}{\tau }, v \right ) & + 
  \delta \left < \frac{ w^{n+1}-w^{n}}{\tau }, v \right >  \\
  & +   \alpha d (w^{\sigma(n)},v) = 0,  
  \quad n = 0,1, ..., N-1,
\end{split}
\end{equation}
\begin{equation}\label{15}
  <w^0,  v> = <\delta^{-\alpha} f , v> ,
  \quad \forall v \in V^h . 
\end{equation} 
where we use the notation
\[
  t^{\sigma(n)} = \sigma t^{n+1} + (1-\sigma) t^{n},
  \quad w^{\sigma(n)} = \sigma w^{n+1} + (1-\sigma) w^{n}.
\]
For $\sigma =0.5$, the difference scheme (\ref{14}), (\ref{15}) approximates the problem 
(\ref{12}), (\ref{13})
with the second order by $\tau$, whereas for other values ​​of $\sigma$, we have only the first order.

The difference scheme (\ref{14}), (\ref{15}) is unconditionally stable under the standard restrictions
(see, e.g., \cite{Samarskii1989,SamarskiiMatusVabischevich2002})  on 
the weight $\sigma \geq 0.5$. To prove this fact, in equation (\ref{14}), we put
\[
 v = \frac{w^{n+1}+w^{n}}{2} + 
 \left (\sigma - \frac{1}{2} \right ) \tau \frac{w^{n+1}-w^{n}}{\tau} .
\] 
For $\sigma \geq 0.5$, we have
\[
 \|w^{n+1}\| \leq \|w^{n}\|,
  \quad n = 0,1, ..., N-1 . 
\]
Using the initial condition (\ref{15}), we get the required a priori estimate for stability
of the scheme (\ref{14}), (\ref{15}) with respect to the initial data:
\begin{equation}\label{16}
 \|w^{n+1}\| \leq \delta^{-\alpha} \|\varphi\|,
  \quad n = 0,1, ..., N-1 . 
\end{equation} 

For $0.5 < \alpha < 1$, we can get an a priori estimate for the difference scheme 
(\ref{14}), (\ref{15}), which is similar to the estimate (\ref{10}) that holds for the differential problem 
(\ref{12}), (\ref{13}). For the sake of simplicity, 
let us consider the scheme (\ref{14}), (\ref{15}) for $\sigma = 0.5$. 
In this case, equation (\ref{14}) may be written as
\begin{equation}\label{17}
  d (t^{n+1/2} w_t, v) + \delta <w_t,v>
  + \alpha d(\widetilde{w},v) = 0,
  \quad n = 0,1, ..., N-1,
\end{equation} 
where
\[
 w_t =  \frac{w^{n+1} - w^{n}}{\tau },
 \quad  \widetilde{w} = \frac{w^{n+1} + w^{n}}{2}.
\] 

Similarly to the problem (\ref{8}), (\ref{9}), in equation (\ref{17}), 
we select
\[
  v = \widetilde{w} + \frac{2\alpha -1}{\alpha} t^{n+1/2} w_t.
\]
This gives
\[
\begin{split}
 \frac{2\alpha -1}{\alpha} t^{n+1/2} \left ( d(t^{n+1/2} w_t, w_t) + \delta <w_t,w_t> \right ) \\
 + \left ( \delta <w_t, \widetilde{w}> + 2\alpha d(t^{n+1/2} w_t,\widetilde{w}\right )
 + \alpha d (\widetilde{w}, \widetilde{w}) = 0 .
\end{split}
\] 
For $0.5 < \alpha < 1$, we have the inequality
\begin{equation}\label{18}
 \delta <w_t,\widetilde{w}> + 2 \alpha S \leq 0, 
\end{equation} 
where
\[
 S =  d(t^{n+1/2} w_t,\widetilde{w}) + \frac{1}{2} d(\widetilde{w}, \widetilde{w}) .
\]
For the first term in (\ref{18}), With the above-mentioned notation, we have
\begin{equation}\label{19}
 \delta <w_t,\widetilde{w}> = \frac{\delta}{2\tau } \left ( \|v^{n+1}\|^2 - \|v^{n+1}\|^2 \right ) . 
\end{equation} 
Direct calculations yield
\[
\begin{split}
 S & = \frac{1}{2\tau } \left ( d( t^{n+1} w^{n+1}, w^{n+1}) - d( t^{n} w^{n}, w^{n}) \right ) \\
 & - \frac{1}{8} d(w^{n+1}, w^{n+1}) + \frac{3}{8} d(w^{n}, w^{n}) + \frac{1}{4} d(w^{n+1}, w^{n}) .
\end{split}
\]
Similarly to the proof for the estimate (\ref{16}), we have
\[
 d(w^{n+1},w^{n+1}) \leq  d(w^{n}, w^{n}) .
\]
By virtue of this, we have
\begin{equation}\label{20}
\begin{split}
 S & \geq \frac{1}{2\tau } \left ( d( t^{n+1} w^{n+1}, w^{n+1}) - d( t^{n} w^{n}, w^{n}) \right ) \\
  & \quad + \frac{1}{8} \left ( d(w^{n+1},w^{n+1}) + 2 d(w^{n+1}, w^{n}) + d(w^{n}, w^{n}) \right ) \\
  & \geq \left ( d( t^{n+1} w^{n+1}, w^{n+1}) - d( t^{n} w^{n}, w^{n}) \right ) .
\end{split} 
\end{equation} 
The substitution of (\ref{19}), (\ref{20}) into (\ref{18}) gives the inequality
\[
 \delta \|w^{n+1}\|^2 + 2\alpha d( t^{n+1} w^{n+1}, w^{n+1}) \leq \delta  \|w^{n}\|^2 + 
 2\alpha d( t^{n} w^{n},w^{n}) .
\] 
Thus, we arrive at the estimate
\begin{equation}\label{21}
 \delta \|w^{n+1}\|^2 + 2\alpha d( t^{n+1} w^{n+1}, w^{n+1}) \leq \delta  \|w^{0}\|^2,
  \quad n = 0,1, ..., N-1 .  
\end{equation}
 
The a priori estimate (\ref{21}) is a discrete analog of the estimate (\ref{6})
for the solution of the Cauchy problem for equation (\ref{8}), (\ref{9}).

Thus, we proved the following basic statement.

\begin{thm}\label{t-2}
For $\sigma \geq 0.5$, the difference scheme (\ref{14}), (\ref{15}) is unconditionally stable
with respect to the initial data.
For $0 < \alpha <1$, the approximate solution satisfies the estimate (\ref{16}), and
for $0.5 < \alpha <1$, we have the estimate (\ref{21}).
\end{thm}

The computational implementation of the difference scheme (\ref{14}), (\ref{15}) involves
the solution of standard elliptic boundary value problems
\[
  d(t^{\sigma(n)}  + 2\alpha \sigma \tau) w^{n+1}, v)  + \delta <w^{n+1}, v> = <\psi^n, v>,
  \quad \forall v \in V^h , 
\] 
for the given $\psi^n$ для $n = 0,1, ..., N-1$. 

\section{Numerical experiments} 

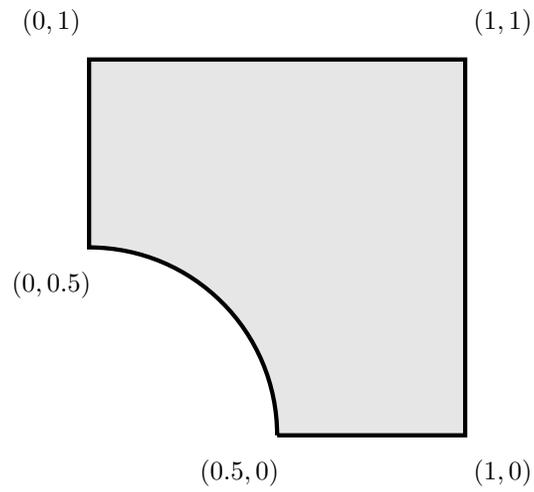
\begin{figure}[htp]
  \begin{center}
  \begin{tikzpicture}[scale = 0.5]
     \draw [ultra thick, fill=gray!20] (5,0) arc [radius=5, start angle=0, end angle= 90] -- (0, 5) -- (0,10) -- (10,10) -- (10,0) -- (5,0);
     \draw(4,-1) node {$(0.5,0)$};
     \draw(11,-1) node {$(1,0)$};
     \draw(-1,4) node {$(0,0.5)$};
     \draw(-1,11) node {$(0,1)$};
     \draw(11,11) node {$(1,1)$};
  \end{tikzpicture}
  \caption{Computational domain $\Omega$}
  \label{f-1}
  \end{center}
\end{figure}

\begin{figure}[!h]
  \begin{center}
    \includegraphics[width=1.\linewidth] {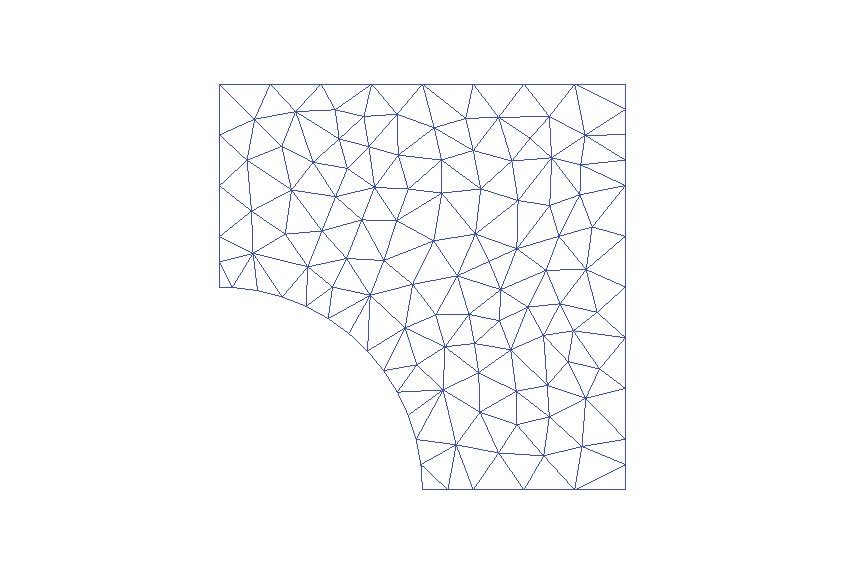}
	\caption{Coarse grid 1: 123 nodes, 202 triangles}
	\label{f-2}
  \end{center}
\end{figure}

\begin{figure}[!h]
  \begin{center}
    \includegraphics[width=1.\linewidth] {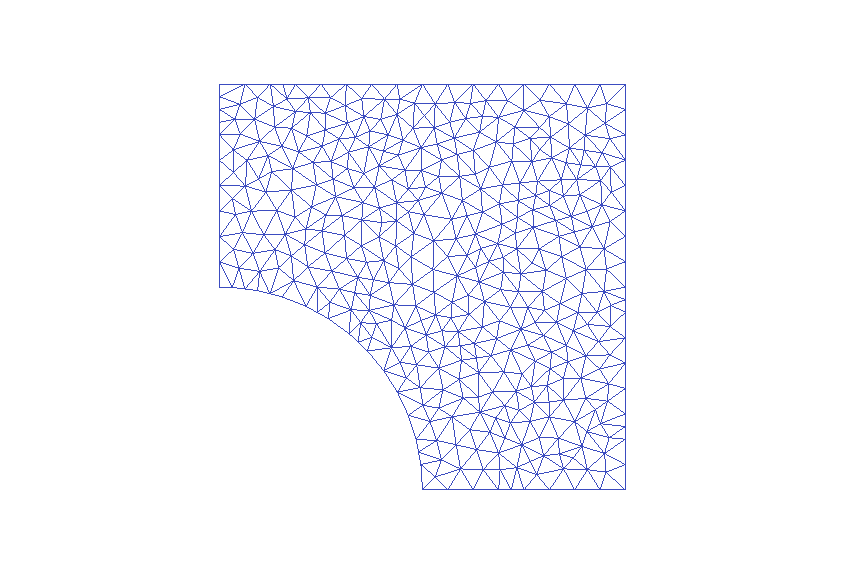}
	\caption{Medium grid 2: 424 nodes, 774 triangles}
	\label{f-3}
  \end{center}
\end{figure}

\begin{figure}[!h]
  \begin{center}
    \includegraphics[width=1.\linewidth] {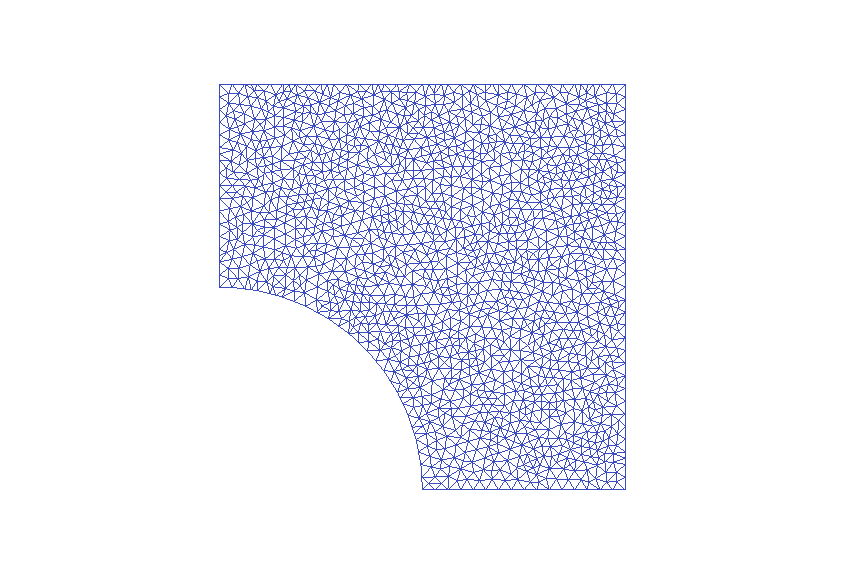}
	\caption{Fine grid 3: 1632 nodes, 3125 triangles}
	\label{f-4}
  \end{center}
\end{figure}

Capabilities of the proposed method are illustrated by solving a two-dimensional model 
problem. The computational domain is a part of the unit square with 
the circular cut; it is shown in Fig.~\ref{f-1}. 
Triangulation is performed to discretize this domain. 
Calculations are performed using  coarse (see Fig.~\ref{f-2}), 
medium (Fig.~\ref{f-3}) and fine (Fig.~\ref{f-4}) grids.

The problem (\ref{4}) for the elliptic operator (\ref{1}), (\ref{2})
is considered with constant coefficients:
\[
 k({\bm x}) = 1, 
 \quad c({\bm x}) = 0,
 \quad \mu ({\bm x}) = \mu. 
\]
Finite element approximations lead to  equation (\ref{5}), (\ref{6}). 
    
To estimate the constant $\delta$, we solve the spectral problem 
\begin{equation}\label{22}
 a(v,v) = \lambda <v,v>,
 \quad v \in V^h,  
\end{equation} 
where $\delta = \lambda_{\min}$.
When choosing  piecewise linear finite elements 
($V^h \subset H^1(\Omega)$), the corresponding values of the constant $\delta$ 
for the above-mentioned computational grids and
$\mu = 1, 10, 100$ are presented in Table~\ref{tbl-1}.

\begin{table}[!h]
\caption{Constant $\delta$}
\begin{center}
\begin{tabular}{|c|c|c|c|} \hline
Сетка    & \multicolumn{3}{c|}{ $\mu$ } \\ \cline{2-4}
    & 1 & 10 & 100  \\ \hline
 1   & 4.09878836985 & 18.2510702649 &  27.1806096252 \\
 2   & 4.09510824780 & 18.1347466168 &  26.7894188262 \\
 3   & 4.09402262307 & 18.1001727374 &  26.6779648687 \\ \hline
\end{tabular}
\end{center}
\label{tbl-1}
\end{table}

\begin{figure}[!h]
  \begin{center}
    \includegraphics[width=0.7\linewidth] {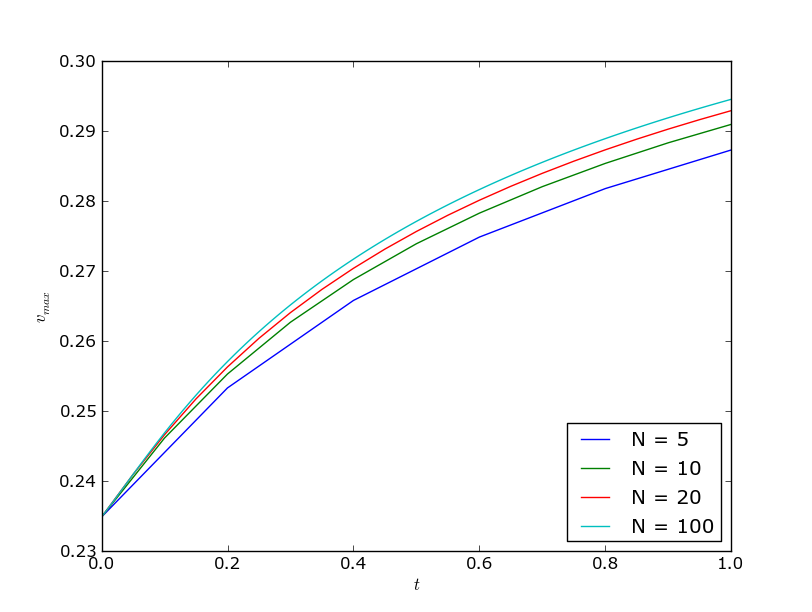}
	\caption{Time-histories of the solution for $\sigma =1$}
	\label{f-5}
  \end{center}
\end{figure}

\begin{figure}[!h]
  \begin{center}
    \includegraphics[width=0.7\linewidth] {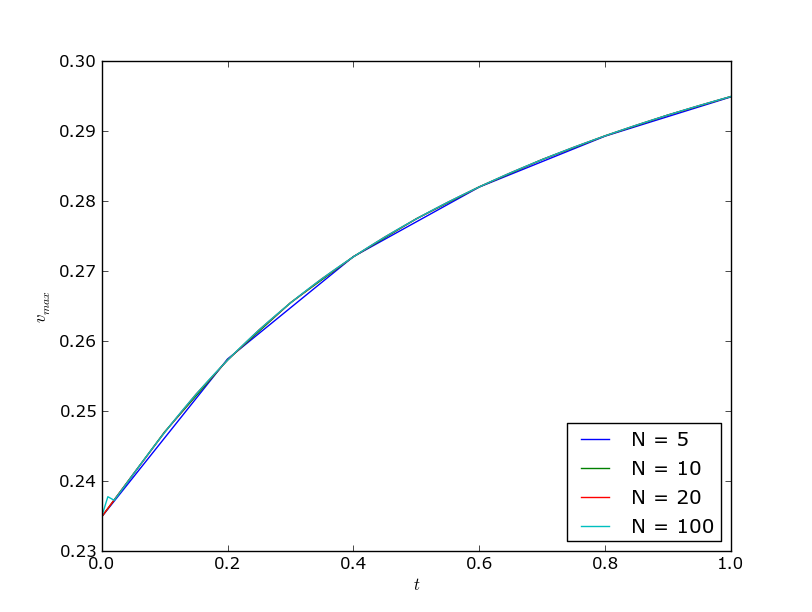}
	\caption{Time-histories of the solution for $\sigma =0.5$}
	\label{f-6}
  \end{center}
\end{figure}

The most interesting fact for this problem is the dependence on the time step.
Figure~\ref{f-5} shows the dependence of the maximum (over the entire computational domain) value
of the approximate solution on the time step.
Calculations are performed on grid 2 using  piecewise linear approximations in space
with $\mu =10$ и $f({\bm x}) = 1$. Unless stated otherwise, we used $\alpha =0.5$.
The convergence of the approximate solution with the first order in time is observed in this figure.
Similar data are presented in Fig.~\ref{f-6} for the symmetric scheme with weights ($\sigma = 0.5$).
Here we see much more rapid convergence and so
we can obtain acceptable by accuracy results using fairly coarse meshes.

\begin{figure}[!h]
  \begin{center}
    \includegraphics[width=1.\linewidth] {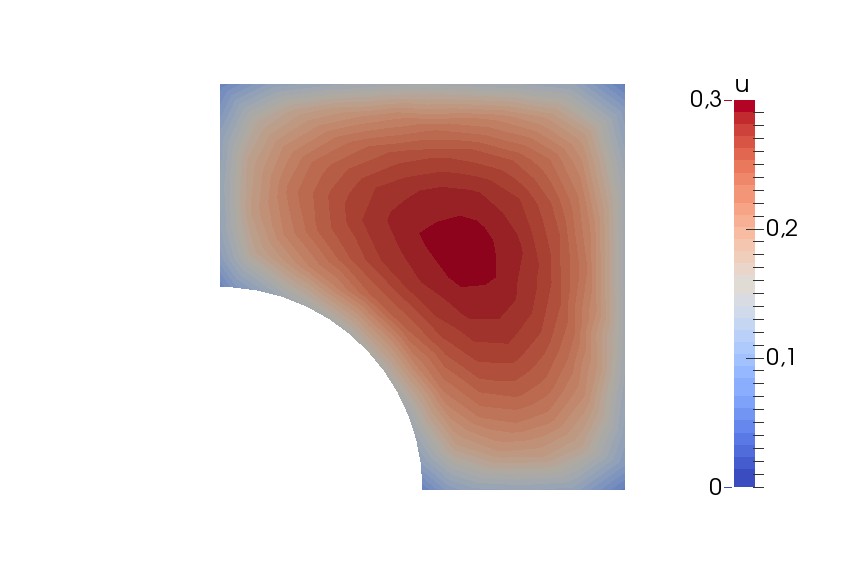}
	\caption{Solution on grid 1 ($v_{max} = 0.294827$)}
	\label{f-7}
  \end{center}
\end{figure}

\begin{figure}[!h]
  \begin{center}
    \includegraphics[width=1.\linewidth] {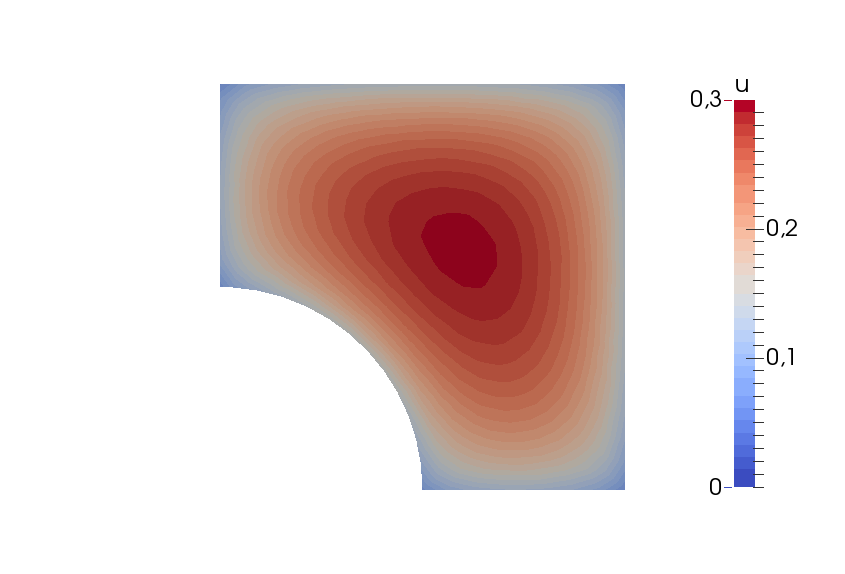}
	\caption{Solution on grid 2 ($v_{max} = 0.294904$)}
	\label{f-8}
  \end{center}
\end{figure}

\begin{figure}[!h]
  \begin{center}
    \includegraphics[width=1.\linewidth] {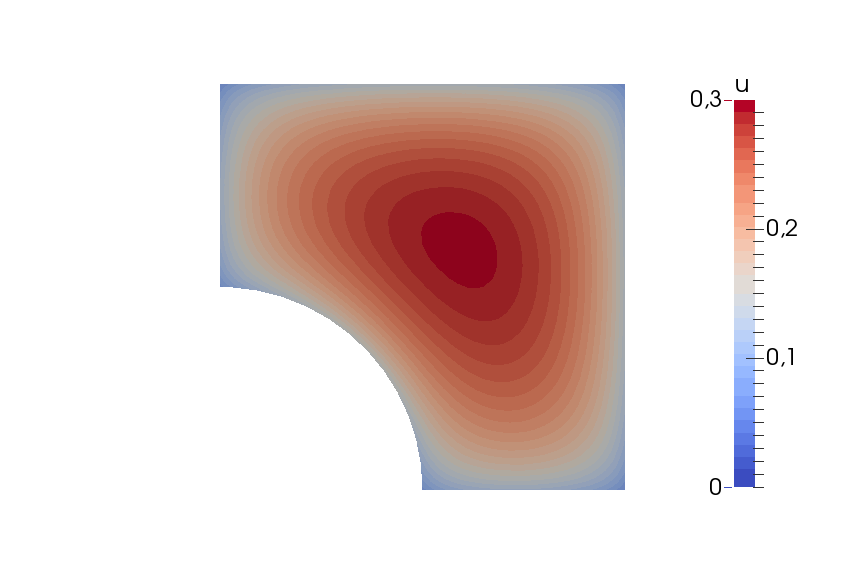}
	\caption{Solution on grid 3 ($v_{max} = 0.29491$)}
	\label{f-9}
  \end{center}
\end{figure}

The dependence on spatial grids for the problem with $\mu =10$ и $f({\bm x}) = 1$  
is depicted in Figs.~\ref{f-7}-\ref{f-9}.
These calculations were performed using the symmetric scheme with the small time step 
that corresponds to $N = 100$.

\begin{figure}[!h]
  \begin{center}
    \includegraphics[width=0.7\linewidth] {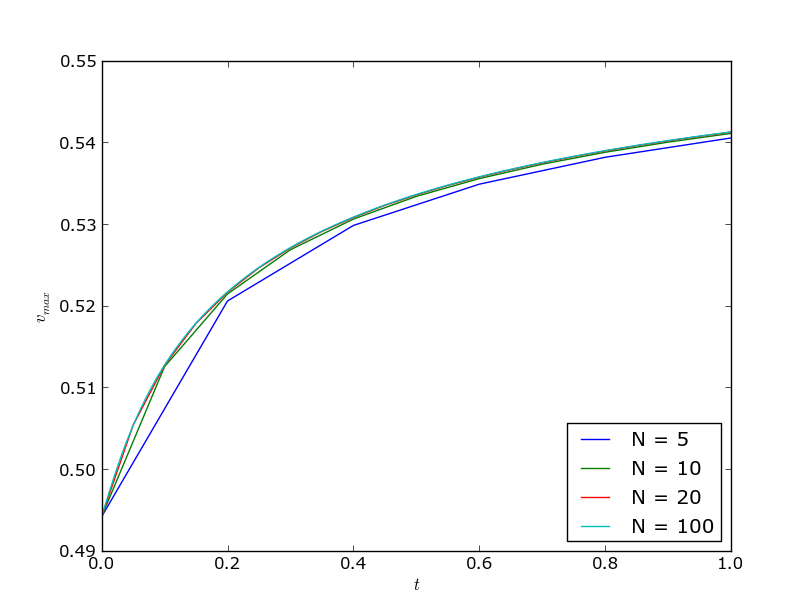}
	\caption{Time-histories of the solution for $\mu  = 1$}
	\label{f-10}
  \end{center}
\end{figure}

\begin{figure}[!h]
  \begin{center}
    \includegraphics[width=0.7\linewidth] {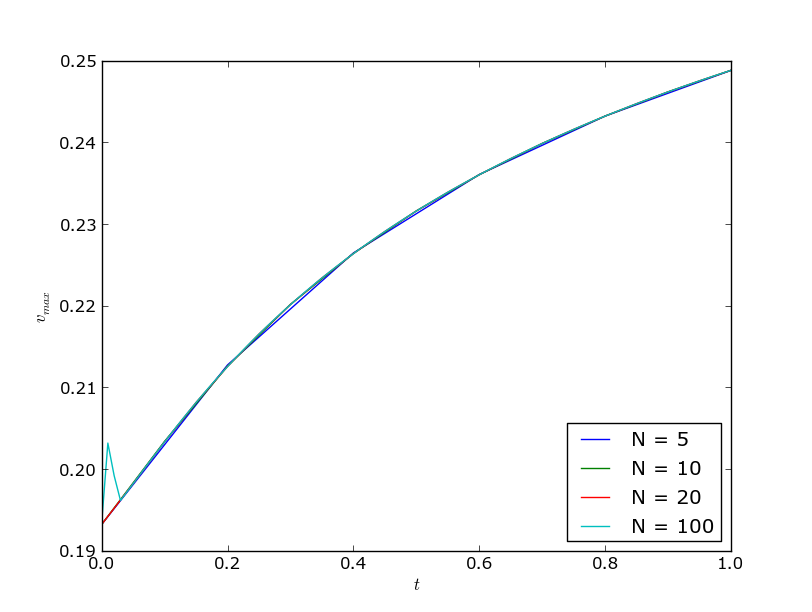}
	\caption{Time-histories of the solution for $\mu  = 100$}
	\label{f-11}
  \end{center}
\end{figure}

The convergence of the symmetric scheme for different values ​​of $\mu$ 
is observed in Figs.~\ref{f-10},\ref{f-11} (see also Fig.~\ref{f-6}).

\begin{figure}[!h]
  \begin{center}
    \includegraphics[width=0.7\linewidth] {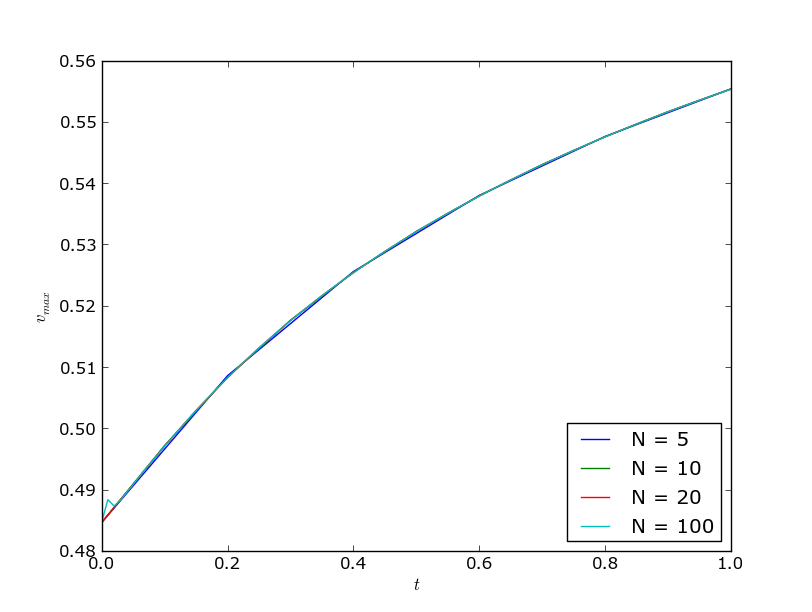}
	\caption{Maximum of the solution for $\alpha  = 0.25$}
	\label{f-12}
  \end{center}
\end{figure}

\begin{figure}[!h]
  \begin{center}
    \includegraphics[width=0.7\linewidth] {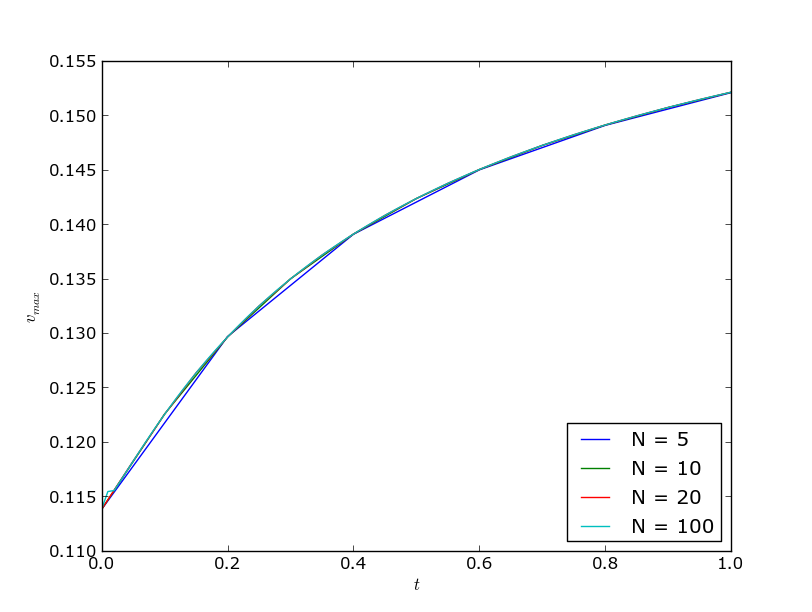}
	\caption{Maximum of the solution for $\alpha   = 0.75$}
	\label{f-13}
  \end{center}
\end{figure}

\begin{figure}[!h]
  \begin{center}
    \includegraphics[width=1.\linewidth] {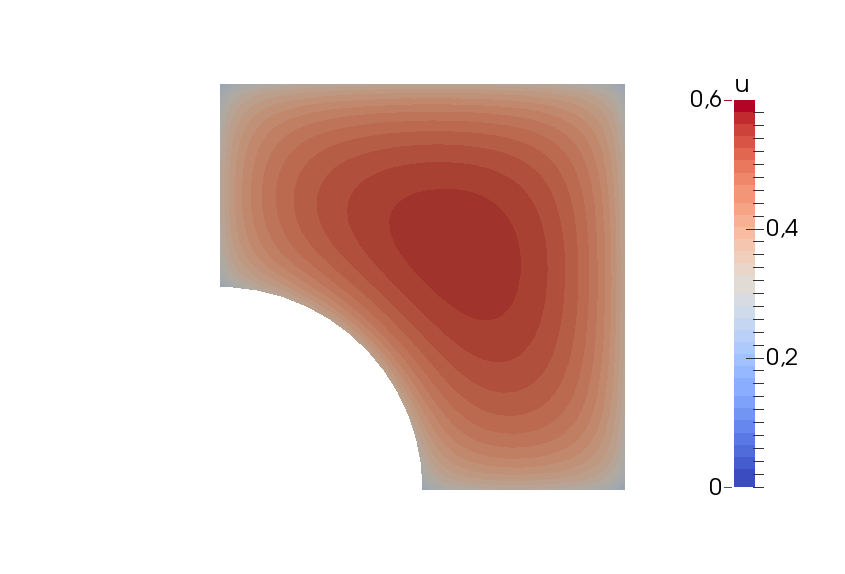}
	\caption{Solution for $\alpha = 0.25$ ($v_{max} = 0.555307$)}
	\label{f-14}
  \end{center}
\end{figure}

\begin{figure}[!h]
  \begin{center}
    \includegraphics[width=1.\linewidth] {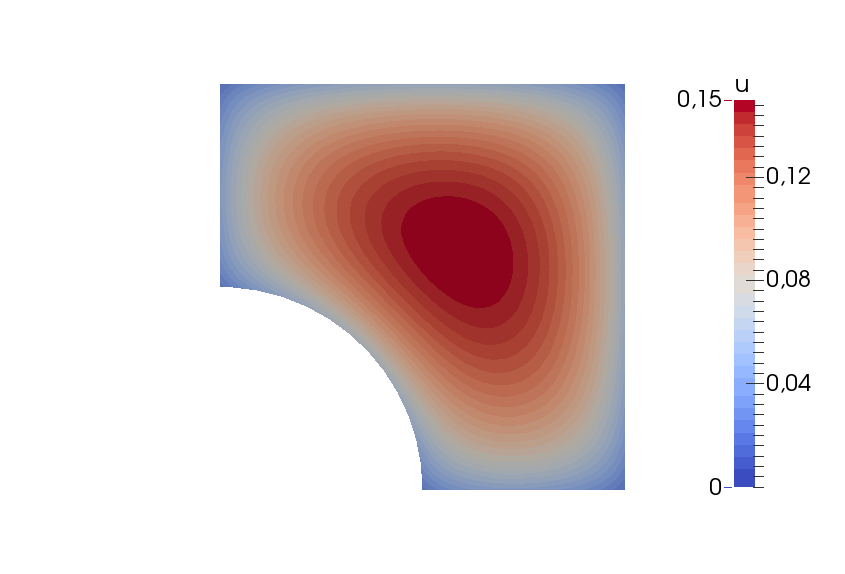}
	\caption{Solution for $\alpha = 0.75$ ($v_{max} = 0.152168$)}
	\label{f-15}
  \end{center}
\end{figure}

It seems reasonable to study the influence of the parameter $\alpha$. 
The convergence of the symmetric scheme for $\alpha = 0.25$ is shown in Fig.~\ref{f-12}, 
similar data for the problem with  $\alpha = 0.75$ are presented in Fig.~\ref{f-13}.
Figures \ref{f-14},~\ref{f-15} demonstrate the solution for various $\alpha$ .

\begin{figure}[!h]
  \begin{center}
    \includegraphics[width=0.7\linewidth] {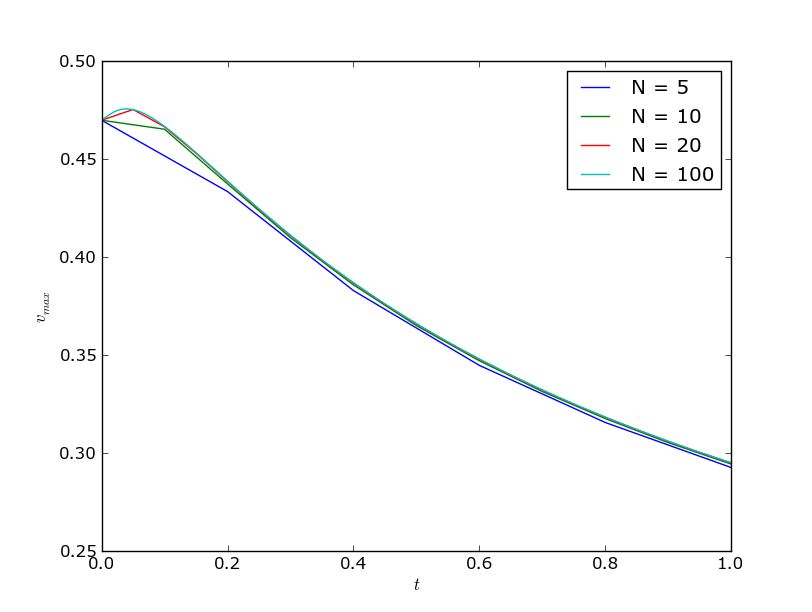}
	\caption{Solution of the problem for $\delta = 0.5 \lambda_{\min}$}
	\label{f-16}
  \end{center}
\end{figure}

\begin{figure}[!h]
  \begin{center}
    \includegraphics[width=0.7\linewidth] {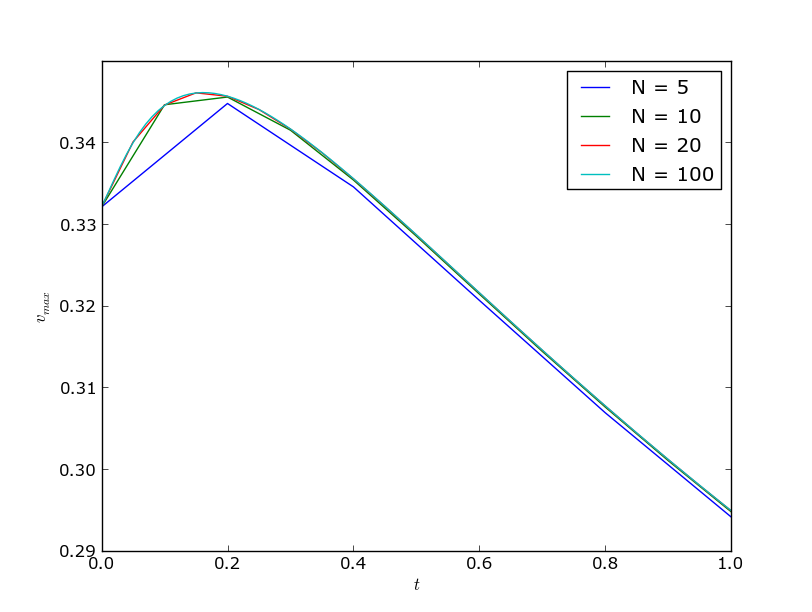}
	\caption{Solution of the problem for $\delta = 0.25 \lambda_{\min}$}
	\label{f-17}
  \end{center}
\end{figure}

Among the main parameters of the problem we can highlight $\delta$.
Features of the computational algorithm (equation (\ref{9})) depend explicitly on this parameter,
although the solution itself is independent of $\delta$.
In these calculations, $\delta$ is evaluated as the principal eigenvalue of
the spectral problem (\ref{22}) ($\delta = \lambda_{\min}$).
But we can also use a rough value for $\delta$.
Figure~\ref{f-16} shows the results for $\delta = 0.5 \lambda_{\min}$, whereas
Fig.~\ref{f-17} corresponds to the case with $\delta = 0.25 \lambda_{\min}$.
We see clearly that the significant inaccuracy in setting $\delta$ does not lead to any significant loss
of efficiency of the computational algorithm.

\clearpage

\section*{Acknowledgements}

This research was supported by RFBR (project 14-01-00785).

\end{document}